\documentclass[12pt]{amsart}
\usepackage[english]{babel}
\usepackage{egothic}
\usepackage[T1]{fontenc}

\usepackage{amsmath}
\usepackage{amsthm}
\usepackage{amssymb}
\usepackage{amsfonts} 
\usepackage{mathrsfs} 
\usepackage{pb-diagram} 
\usepackage[all]{xy}
\usepackage{mathrsfs}
\usepackage{enumerate}
\usepackage{physics}
\usepackage[notcite, final, notref]{showkeys}
\usepackage{dsfont}
\usepackage[dvips]{color}
\newtheorem{theorem}{Theorem}[section]

\newtheorem{proposition}[theorem]{Proposition}

\newtheorem{definition}[theorem]{Definition}

\theoremstyle{definition}
\newtheorem{remark}[theorem]{Remark}

\newtheorem{example}[theorem]{Example}

\newcommand{\C}{{\mathbb C}}
\newcommand{\R}{{\mathbb R}}

\usepackage{xcolor}

\begin{document}

\title[New Fueter Variables, Global Operator, Quaternionic Case]{New Fueter Variables Associated to the Global Operator in the Quaternionic Case}

\author[D. Alpay]{Daniel Alpay}
\address{(DA) Schmid College of Science and Technology \\
Chapman University\\
One University Drive
Orange, California 92866\\
USA}
\email{alpay@chapman.edu}

\author[K. Diki]{Kamal Diki}
\address{(KD) Schmid College of Science and Technology \\
Chapman University\\
One University Drive
Orange, California 92866\\
USA}
\email{diki@chapman.edu}

\author[M. Vajiac]{Mihaela Vajiac}
\address{(MV) Schmid College of Science and Technology \\
Chapman University\\
One University Drive
Orange, California 92866\\
USA}
\email{mbvajiac@chapman.edu}

\thanks{D. Alpay thanks the Foster G. and Mary McGaw Professorship in
  Mathematical Sciences, which supported his research. D. Alpay and M. Vajiac thank Chapman University for its Faculty Opportunity Fund that helped fund this research. K. Diki thanks Prof. I. Sabadini and Prof. F. Colombo for introducing him to the theory of slice hyperholomorphic functions using the global operator approach and the Grand Challenges Initiative Program at Chapman University for his postdoctoral position that enabled this work.}
 
\begin{abstract}
The purpose of this paper is to develop a new theory of three non-commuting quaternionic variables and its related Schur analysis theory for a modified version of the quaternionic global operator.
\end{abstract}

\maketitle

\noindent {\em }
\date{today}

{\em Please note that no datasets were generated or analyzed during the current study.}

\section{Introduction}
In 2007 a new theory of regular functions of a quaternionic variable was introduced in \cite{GS2007}. These are the so-called {\em slice hyperholomorphic} functions on quaternions. For more details on this topic we refer the reader to the books \cite{CSS2016, CSS2011, GSS2022} and the references therein. This theory generated several interesting applications in different areas of mathematics and physics due to the discovery of the notion of the S-spectrum for quaternionic operators. In particular, this function theory allowed to develop the quaternionic counterpart of operator theory, Schur analysis, and quantum mechanics, see \cite{ACS2016Schur, CG2019, CGK2018, CSS2011}. A fundamental technique that we use in this paper is based on the research developed in \cite{CGCS2013} where the authors discovered a new approach to the theory of slice hyperholomorphic functions using a special operator with non-constant coefficients called the {\em quaternionic global operator}.\\
The purpose of this paper is to develop a theory of three non-commuting quaternionic variables and a related Schur analysis for the quaternionic global operator. We first set the framework for our work, and define the skew field of quaternions (denoted by $\mathbb H$) as the space of elements of the form:
\begin{equation}
  \label{q123}
  q=x_0+x_1{\mathbf e_1}+x_2{\mathbf e_2}+x_3{\mathbf e_3},
  \end{equation}
  where $x_0,x_1,x_2,x_3$ are real numbers and
  $\mathbf e_1,\mathbf e_2,\mathbf e_3$ satisfy the Cayley multiplication table. A number of the form $q$ in \eqref{q123} is called a quaternion, and
\[
  \vec{q}=x_1\mathbf e_1+x_2\mathbf e_2+x_3\mathbf e_3
\]
is called its vector part (more information on quaternions is provided in Section \ref{sec-2}). 

As in Definition~\ref{G_q_def}, the global operator with non-constant coefficients, which was introduced for the first time in \cite{CGCS2013}, was initially written as:
\begin{equation}
\label{G_q}
\displaystyle G_q(f):=|\vec{q}\,|^2\partial_{x_0}f(q)+\vec{q}\,\sum_{l=1}^3x_l\partial_{x_l}f(q).
\end{equation}

The operator $G_q$ acts on functions of a quaternionic variable which are differentiable with respect to the real variables and it was first used to develop a Cauchy formula in the case of slice hyperholomorphic functions, see~\cite{CGCS2013}. For a distributional approach to the Cauchy problem in this context see~\cite{CS2014}.
In particular, in~\cite{CGCS2013}, it was proved that slice-hyperholomorphic functions are strictly included in the kernel of the global operator $G_q$. 
However, the global operator we work with, denoted by $V_q$, is the normalized form of $G_q$, $\displaystyle V_q=\frac{G_q}{|\vec{q}\,|^{2}}$, namely:
\begin{equation}
  \label{vqbelow}
  V_q=\frac{\partial}{\partial x_0}-\frac{1}{\vec{q}}\sum_{u=1}^3x_u\frac{\partial}{\partial x_u}.
\end{equation}

The global operator $V_q$ was used  in~\cite{ADS2021}  to develop the Fueter mapping theorem for poly slice monogenic functions.  In~\cite{ACPSS} the global operator was used to study infinite-order differential operators acting on entire hyperholomorphic functions and it was used also in \cite{GP2014} to study global differential equations for slice regular functions. Recent extensions and new results on the global operator and related topics can be found in papers \cite{GCGC2021, GCGC2021-2}.

\smallskip

\begin{definition} A quaternionic-valued function $f$ which is
at least $C^1$ in the real variables   in a domain $\Omega_R\subset\mathbb R^4$ is called
  $V_q$-regular if $V_qf=0$.
\end{definition}

In this paper we will focus on a family of $V_q$-regular functions which are real-analytic and we develop a Cauchy-Kovalevskaia product (denoted CK-product) for this space.  We consider the version of Gleason problem associated to $V_q$, which will allow us to introduce a new type of Fueter-like
variables and develop a Schur analysis and rational functions for this operator. For consistency we now introduce the notations that will be used throughout the paper. 

\begin{definition}
\label{V_q-Fueter_var}
Given $x=(x_0,x_1,x_2,x_3)\in\mathbb{R}^4$ such that $x_u\neq 0$ for some $u=1,2,3$. We define the $V_q$-Fueter variables:
\begin{equation}
\mu_u(x)=x_u\left(1+\frac{x_0}{\vec{q}}\right),\quad u=1,2,3,
\end{equation}
and, for $\alpha=(\alpha_1,\alpha_2,\alpha_3)\in\mathbb{N}_0^3$, we write, with a slight abuse of notation:
$$
x^{\alpha}=x_1^{\alpha_1}x_2^{\alpha_2}x_3^{\alpha_3},
$$
where $x^{\alpha}$ does not depend on $x_0$.
\end{definition}

The functions $\mu_u$, $u=1,2,3$, are the counterpart of the classical Fueter variables (see Definition~\ref{classical-F} for the latter) in the present setting.  These new variables allow us  to define
and study a counterpart of Schur analysis in this
setting. 
In Theorem~\ref{KerV_q} we prove a more general form of the following:

\begin{theorem}
The functions
\begin{equation}
\mu_u(x)=x_u\left(1+\frac{x_0}{\vec{q}}\right),\quad u=1,2,3,
\end{equation}
are $V_q$-regular on $\mathbb{H}^*:= \lbrace q=x_0+\vec{q}, \quad \vec{q}\neq 0\rbrace$. 
\end{theorem}

Throughout the paper we will exchange the quadruple $x=(x_0,x_1,x_2,x_3)$ with $x_u\neq 0$ for some $u=1,2,3$ with the quaternionic variable $q=x_0+\vec{q}\in \mathbb{H}^*$ so that we can consider $\mu_u$ as functions of both. For a given $x\in\mathbb{H}^*$, $\mu_u(x)$ commute pairwise and we recover Fueter-like monomials as in Definition~\ref{mu-alpha}: 
\begin{equation}
  \mu^\alpha(x)=  \mu_1^{\alpha_1}(x)  \mu_2^{\alpha_2}(x)  \mu_3^{\alpha_3}(x),
\end{equation}
where $\alpha=(\alpha_1,\alpha_2,\alpha_3)\in\mathbb N_0^3$.

The functions $\mu^\alpha$ are $V_q$-regular on $\mathbb{H}^*$ (see Section~\ref{Fueter-def}) and we will see that they can be rewritten as
\begin{equation}
\mu^\alpha(x)=x^\alpha\left(1+\frac{x_0}{\vec{q}}\right)^{|\alpha|}.
\end{equation}

As explained in the following remarks, and throughout
the paper there are important and fundamental differences between the cases considered in
\cite{MR3819695,MR2124899,MR3795166} (and in particular the
case of Fueter variables) and the present setting.

\begin{remark}
  In the setting of Fueter variables, the Cauchy-Dirac operator and the operators of partial differentiation commute for smooth functions.
    This is not the case for the operator $V_q$, which does not commute with differentiations with respect to $x_u$, $u=1,2,3$. So, for instance,
    $\displaystyle\frac{\partial f}{\partial x_1}$ need not be $V_q$-regular when $f$ is $V_q$-regular. Since $V_q$ and
    $\displaystyle\frac{\partial f}{\partial x_0}$ commute,  $\displaystyle\frac{\partial f}{\partial x_0}$ is $V_q$-regular when $f$ is $V_q$-regular.
    As an illustration the case $f(x)=\mu_1(x)$ is considered in Example \ref{totoche}.
    The function
    $\displaystyle\frac{\partial\mu_1}{\partial x_2}$ is not $V_q$-regular, while the function
$\displaystyle\frac{\partial\mu_1}{\partial x_0}$ is $V_q$ -regular but cannot be written as a convergent Fueter-like series of the $\mu^{\alpha}$ monomials.
  \end{remark}

To present these ideas we adopt the following structure: in Section~\ref{prelim} we review basic notations and definitions of quaternions, slice hyperholomorphic functions, and the quaternionic global operator with nonconstant coefficients. In Section~\ref{Fueter-def} we introduce the Fueter-like variables and prove a Gleason type theorem in this setting. Section~\ref{CK-Vq} is devoted to the study of a Cauchy-Kovalevskaia type extension and CK-type product corresponding to these Fueter-like variables. In Section~\ref {Vq-Kernel-Omega} and~\ref{Vq-Arveson} we study some examples of reproducing kernel Hilbert spaces generated by the Fueter-like variables including a counterpart of Averson space and Schur multipliers. Finally, in Section~\ref{Vq-Blaschke} we obtain Blaschke factors in this framework followed by building a rational function theory in Section~\ref{V_q-Rational}.
We conclude the paper with a general description of the next steps we will undertake using these new Fueter-like variables.

\section{Preliminary results}
\setcounter{equation}{0}
\label{prelim}
We recall that the non-commutative field of quaternions is defined and denoted by 
$$\mathbb{H}=\lbrace{q=x_0+x_1i+x_2j+x_3k : \ x_0,x_1,x_2,x_3\in\R}\rbrace,$$ where the imaginary units satisfy the multiplication rules $$i^2=j^2=k^2=-1\quad \text{and}\quad ij=-ji=k, jk=-kj=i, ki=-ik=j.$$
On $\mathbb{H}$ the conjugate and the modulus of $q$ are defined respectively by
$$\overline{q}=x_0-\vec{q}\,,  \quad \vec{q}\,=x_1i+x_2j+x_3k$$
and $$\vert{q}\vert=\sqrt{q\overline{q}}=\sqrt{x_0^2+x_1^2+x_2^2+x_3^2}.$$

Throughout the paper we use the notations ${\mathbf e}_0=1, {\mathbf e}_1=i, {\mathbf e}_2=j$ and ${\mathbf e}_3=k$ interchangeably for the imaginary units. It is important to note that the conjugation satisfy the property $\overline{ pq }= \overline{q}\, \overline{p}$ for any $p,q\in \mathbb{H}$.

An important tool in the quaternionic case is the symmetric product of $n$ quaternionic numbers $q_1, q_2,\dots, q_n$:
\begin{equation}
      \label{sumsum}
    q_1\times q_2\times\cdots\times q_n=\frac{1}{n!}\sum_{\sigma\in S_n}q_{\sigma(1)}q_{\sigma(2)}\cdots q_{\sigma(n)},
\end{equation}
where the sum is over the set $S_n$ of all permutations on $n$ indices.

In the classical case, Fueter~\cite{fueter-32, MR1509533} used the following variables to describe the kernel of the Cauchy-Fueter operator:
\begin{definition}
\label{classical-F}
The classical Fueter variables are:  $\zeta_l= x_l-x_0{\mathbf e}_l$ and we denote by $\displaystyle \zeta^{\alpha}$ the symmetric product:
$\zeta_1^{\times \alpha_1}\times \zeta_2^{\times \alpha_2}\times \zeta_3^{\times \alpha_3}$, where $\alpha=(\alpha_1, \alpha_2, \alpha_3)$.
\end{definition}
These variables are Fueter-regular with respect to the classical Cauchy-Fueter operator.\\
We now define the setting for slice hyperholomorphic functions. The unit sphere $$\lbrace{\vec{q}=x_1i+x_2j+x_3k : \text{ } x_1^2+x_2^2+x_3^2=1}\rbrace$$ coincides with the set of all  imaginary units given by $$\mathbb{S}=\lbrace{q\in{\mathbb{H}} : q^2=-1}\rbrace.$$

\begin{definition}
Let $f: \Omega \longrightarrow \mathbb{H}$ be a  $\mathcal{C}^1$ function on a given domain $\Omega\subset \mathbb{H}$. Then, $f$ is said to be (left) slice hyperholomorphic function if, for every $I\in \mathbb{S}$, the restriction $f_I$ to $\C_{I}=\R+I\R$, with variable $q=x+Iy$, is holomorphic on $\Omega_I := \Omega \cap \C_I$, that is it has continuous partial derivatives with respect to $x$ and $y$ and the function
$\overline{\partial_I} f : \Omega_I \longrightarrow \mathbb{H}$ defined by
$$
\overline{\partial_I} f(x+Iy):=
\dfrac{1}{2}\left(\frac{\partial }{\partial x}+I\frac{\partial }{\partial y}\right)f_I(x+yI)
$$
vanishes identically on $\Omega_I$. The set of slice hyperholomorphic functions will be denoted by $\mathcal{SR}(\Omega)$.
\end{definition}
The right quaternion vector space $\mathcal{SR}(\Omega)$ is endowed with the natural topology of uniform convergence on compact sets. The characterization of such functions on a ball centered at the origin is given by:

\begin{theorem}[Series expansion \cite{GS2007}]
An $\mathbb{H}$-valued function $f$ is slice hyperholomorphic on $B(0,R)$ if and only if it has a series expansion of the form:
$$f(q)=\sum_{n=0}^{+\infty} q^na_n$$
converging on $B(0,R)=\{q\in\mathbb{H};\mid q\mid<R\}$.
\end{theorem}
Another interesting approach to define slice hyperholomorphic functions is to consider them as solutions of a special global operator with non constant coefficients that was introduced and studied in \cite{CGCS2013, CS2014, GP2014}. This leads to the following definition
\begin{definition}
\label{G_q_def}
Let $\Omega$ be an open set in $\mathbb{H}$ and $f:\Omega\longrightarrow \mathbb{H}$ a function of class  $\mathcal{C}^1$. We define the global operator $G_q(f)$ by
$$\displaystyle G_q(f):=|\vec{q}\,|^2\partial_{x_0}f(q)+\vec{q}\,\sum_{l=1}^3x_l\partial_{x_l}f(q),$$
for any $q=x_0+\vec{q}\,\in\Omega$.
\end{definition}
It was proved in \cite{CGCS2013} that any slice hyperholomorphic function is in the kernel of $G_q$ on axially symmetric slice domains. We briefly recall the definition of such a domain:
\begin{definition}

A domain $\Omega\subset \mathbb{H}$ is said to be a slice domain (or just $s$-domain) if  $\Omega\cap{\mathbb{R}}$ is nonempty and for all $I\in{\mathbb{S}}$, the set $\Omega_I:=\Omega\cap{\mathbb{C}_I}$ is a domain of the complex plane $\mathbb{C}$.
If moreover, for every $q=x+Iy\in{\Omega}$, the whole sphere $x+y\mathbb{S}:=\lbrace{x+Jy; \, J\in{\mathbb{S}}}\rbrace$
is contained in $\Omega$, we say that  $\Omega$ is an axially symmetric slice domain.
\end{definition}
 There are many other interesting properties of the global operator $G_q$ that were studied in the literature, in particular in~\cite{CGCS2012}. We recall some here:
\begin{proposition}\label{GPRO}
Let $\Omega$ be an open set in $\mathbb{H}$ and $f,g:\Omega\longrightarrow \mathbb{H}$ two functions of class $\mathcal{C}^1$. Then, for $q=x_0+\vec{q}\,\in\Omega$ we have
\begin{enumerate}
\item $\displaystyle G(fg)=G(f)g+fG(g)+(\vec{q}\,f-f\vec{q}\,)\sum_{l=1}^3x_l\partial_{x_l}g$. \newline \text{In particular, it holds:}

\item $G(f\lambda+g)=G(f)\lambda+G(g),  \forall \lambda \in \mathbb{H}$.
\item $G(x_0f)=|\vec{q}\,|^2f+x_0G(f)$ and $G(\vec{q}\,f)=-|\vec{q}\,|^2f+\vec{q}\,G(f).$
\item $G(q^kf)=q^kG(f), \forall k\in\mathbb{N}.$

\end{enumerate}
\end{proposition}

  \section{Fueter-type variables}
  \setcounter{equation}{0}
  \label{Fueter-def}

In this section we introduce the Gleason setting that will yield the new $V_q$-Fueter variables in the quaternionic case, in a natural way. 
Through this process, we find these variables which correspond to the operator $V_q$ when we use the same technique the authors applied in the real ternary case~\cite{MR3795166}, in the split quaternionic case~\cite{MR3819695}, and in the regular Fueter
quaternionic case~\cite{MR2124899}. For completion, in the case of Grassmann variables, see also~\cite{zbMATH07202564}. \smallskip

It is also worth pointing out that our work here is the first that finds a counterpart of the classical Fueter variables through the treatment of this  case of a differential operator with non-constant coefficients. All the applications mentioned before were using this method for operators with constant coefficients.

The strategy involves a clever application of the chain rule, and, for a $V_q$-regular function $f$, we compute:
\[
\int_a^b\frac{\rm d}{\rm dt}f(tx)dt.
\]
In the previous examples of ternary and split quaternions, we can take $a=0$ and $b=1$. 
In the present case $a$ cannot be chosen to be $0$ because of the singularity
at the origin.

We will see that we can apply the same technique to obtain the Fueter-like variables $\mu_u$ in the proof of the following theorem:
\begin{theorem}
Let $\Omega$ be an open domain of $\mathbb{H}^*$.
For a $V_q$-regular function $f\in C^1(\Omega)$, let $a,b\in\Omega$ such that $[a,b]=\{(1-t)a+tb,\, 0\le t\le 1\} \subset\Omega$.
  Then
  \begin{equation}
f(b)-f(a)=\sum_{u=1}^3\mu_u(b-a)R_u^{ab}f,
\end{equation}
where
\begin{equation}
  R_u^{ab}f =\int_a^b\frac{\partial f}{\partial x_u}(a+t(b-a))dt.
  \end{equation}
\end{theorem}

\begin{proof} The chain rule gives
  \[
    \begin{split}
      \frac{\rm d}{\rm dt}f(a+tx)&=x_0\frac{\partial f}{\partial x_0}(a+tx)+\sum_{u=1}^3x_u\frac{\partial f}{\partial x_u}(a+tx),\\
      &=\frac{x_0}{\vec{q}}\sum_{u=1}^3x_u\frac{\partial f}{\partial x_u}(a+tx)+\sum_{u=1}^3x_u\frac{\partial f}{\partial x_u}
      (a+tx),\\
      &=\sum_{u=1}^3\mu_u(x)\frac{\partial f}{\partial x_u}(a+tx),
    \end{split}
    \]
    and the result follows by integrating back and setting $x=b-a$.
   The reader should note that in the second step of the equality we use the fact that $f$ is in the kernel of $V_q$.
  \end{proof}

\begin{remark}{\rm
 Note that
 \begin{equation}
   \mu_u(b-a)\not=\mu_u(b)-\mu_u(a).
 \end{equation}
 Furthermore, the operators $R_u^{ab}$ do not commute with $V_q$.  }
\end{remark}

For $\alpha=(\alpha_1,\alpha_2,\alpha_3)\in\mathbb N_0^3$ we use the multi-index notation and set
\begin{equation}
  x^\alpha=x_1^{\alpha_1}x_2^{\alpha_2}x_3^{\alpha_3}\quad{\rm and}\quad |\alpha|=\alpha_1+\alpha_2+\alpha_3.
\end{equation}

\begin{definition}
\label{mu-alpha}
For the Fueter variables $\mu_u$ and $\alpha=(\alpha_1,\alpha_2,\alpha_3)\in\mathbb N_0^3$, we define the product:
\begin{equation}
  \mu^\alpha(x)=  \mu_1^{\alpha_1}(x)  \mu_2^{\alpha_2}(x)  \mu_3^{\alpha_3}(x),
\end{equation}
for every $x\in\mathbb{H}^*$
\end{definition}

\begin{remark} Since the  $V_q$-Fueter variables $\mu_u$ commute, in the above definition we do not need to use the symmetric product and, in fact, we can rewrite the product as:
\begin{equation}
\mu^\alpha(x)=x^\alpha\left(1+\frac{x_0}{\vec{q}}\right)^{|\alpha|}.
\end{equation}
\end{remark}

\begin{proposition}
 Let $\Omega$ be an open domain in $\mathbb{H}^*$. For every $x\in\Omega\subset \mathbb{H}^*$ it holds that
  \begin{equation}
    |\mu^\alpha|^2=|x|^{2\alpha}\left(1+\frac{x_0^2}{x_1^2+x_2^2+x_3^2}\right)^{|\alpha|},
    \end{equation}
    and, in particular, we have:
    \begin{equation}
      \label{ineq987}
|\mu_u(x)|^2\leq |\zeta_u(x)|^2,\quad u=1,2,3.
      \end{equation}
  \end{proposition}

  \begin{proof}
This follows from the fact that $\frac{x_0}{\vec{q}}$ has no real part on $\Omega$ and:
    \begin{equation}
|\mu_u(x)|^2=x_u^2+\frac{x_u^2x_0^2}{x_1^2+x_2^2+x_3^2}\le x_u^2+x_0^2=|\zeta_u(x)|^2.
      \end{equation}
    \end{proof}
    
We can now prove that the $V_q$-Fueter products $\mu^\alpha$ are in the kernel of $V_q$ on $\mathbb{H}^*$:
\begin{theorem} 
\label{KerV_q}

It holds that $\mu^\alpha$ are in the kernel of the operator $V_q$ on any open domain $\Omega\subset \mathbb{H}^*$. Moreover, we have:
  \begin{equation}
V_q\mu^\alpha(x)=0,
    \end{equation}
    for every $x\in\Omega$.
  \end{theorem}
  
  \begin{proof} We divide the verification into a number of steps. \\

    STEP 1: {\sl It holds that
      \begin{equation}
        \label{kobe}
\frac{\partial }{\partial x_u}\frac{1}{\vec{q}}=\frac{\mathbf e_u}{\vec{q\,}^2}+\frac{2x_u}{\vec{q\,}^3},\quad u=1,2,3.
\end{equation}
}

Indeed, we have
\[
\frac{1}{\vec{q}}=\frac{\vec{q}}{\vec{q\,}^2}=-\frac{\vec{q}}{x_1^2+x_2^2+x_3^2}.
  \]
  Hence
  \[
\begin{split}
  \frac{\partial }{\partial x_u}\frac{1}{\vec{q}}&=\frac{-\mathbf e_u(x_1^2+x_2^2+x_3^2)
    +2x_u\vec{q}}{(x_1^2+x_2^2+x_3^2)^2}\\
  &=\frac{\mathbf e_u\vec{q}^{\,2}}{\vec{q}^{\,4}}+\frac{2x_u\vec{q}}{\vec{q}^{\,3}}
\end{split}
\]
and hence the result.\\


STEP 2: {\sl It holds that
  \begin{equation}
    \label{osaka123}
    \frac{\partial}{\partial x_u}\left(1+\frac{x_0}{\vec{q}}\right)^{|\alpha|}=
    \sum_{\substack{t,s\in\mathbb N_0\\t+s=|\alpha|}}
    \left(1+\frac{x_0}{\vec{q}}\right)^t\left(\frac{x_0\mathbf e_u}{\vec{q\,}^2}+\frac{2x_ux_0}{\vec{q\,}^3}\right)
    \left(1+\frac{x_0}{\vec{q}}\right)^s,
  \end{equation}
  where $u=1,2,3.$}\smallskip

\eqref{osaka123} is a direct consequence of \eqref{kobe} and of the formula for the derivative of $f^n$ when $f$ is a
matrix-valued function (and in particular quaternionic valued) of (say) a real variable $w$:
\begin{equation}
  \label{rtyu}
\frac{{\rm d} f^n}{{\rm d}w}=\sum_{\substack{t,s\in\mathbb N_0\\t+s=|\alpha|}}f^tf^\prime f^s.
\end{equation}

STEP 3: {\sl We have
  \begin{equation}
    \label{osaka12345}
    \frac{\partial}{\partial x_0}\left(1+\frac{x_0}{\vec{q}}\right)^{|\alpha|}=
    \frac{|\alpha|}{\vec{q}}    \left(1+\frac{x_0}{\vec{q}}\right)^{|\alpha|-1}.
  \end{equation}
}
This is because $1+\frac{x_0}{\vec{q}}$ commutes with its derivative with respect to $x_0$ , and formula \eqref{rtyu}
reduces then to the classical formula.\\

STEP 4: {\sl We now calculate
  \[
\frac{1}{\vec{q}}
    \sum_{u=1}^3x_u\frac{\partial}{\partial x_u}\mu^\alpha\, .
    \]
}

We have:
\[
  \begin{split}
    \frac{1}{\vec{q}}
    \sum_{u=1}^3x_u\frac{\partial}{\partial x_u}x^\alpha\left(1+\frac{x_0}{\vec{q}}\right)^{|\alpha|}&=\\
    &\hspace{-4cm}=
\frac{1}{\vec{q}}\left[\sum_{u=1}^3x_u\left(\alpha_ux^{\alpha-e_u}\left(1+\frac{x_0}{\vec{q}}\right)^{|\alpha|}
+x^\alpha\sum_{\substack{t,s\in\mathbb N_0\\t+s=|\alpha|-1}}
\left(1+\frac{x_0}{\vec{q}}\right)^t\left(\frac{x_0\mathbf e_u}{\vec{q\,}^2}+\frac{2x_ux_0}{\vec{q\,}^3}\right)
\left(1+\frac{x_0}{\vec{q}}\right)^s\right)\right]\\
&\hspace{-4cm}=\frac{1}{\vec{q}}\left[|\alpha|x^\alpha\left(1+\frac{x_0}{\vec{q}}\right)^{|\alpha|}
+x^\alpha\sum_{\substack{t,s\in\mathbb N_0\\t+s=|\alpha|-1}}
\left(1+\frac{x_0}{\vec{q}}\right)^t\left(\sum_{u=1}^3\frac{x_ux_0\mathbf e_u}{\vec{q\,}^2}+\frac{2x_u^2x_0}
  {\vec{q\,}^3}\right)
\left(1+\frac{x_0}{\vec{q}}\right)^s\right]\\
&\hspace{-4cm}=\frac{x^\alpha}{\vec{q}}\left[|\alpha|\left(1+\frac{x_0}{\vec{q}}\right)^{|\alpha|}
+\sum_{\substack{t,s\in\mathbb N_0\\t+s=|\alpha|-1}}
\left(1+\frac{x_0}{\vec{q}}\right)^t\left(\frac{x_0\vec{q}}{\vec{q\,}^2}-\frac{2\vec{q\,}^2x_0}
  {\vec{q\,}^3}\right)
\left(1+\frac{x_0}{\vec{q}}\right)^s\right]\\
&\hspace{-4cm}=\frac{x^\alpha}{\vec{q}}\left[|\alpha|\left(1+\frac{x_0}{\vec{q}}\right)^{|\alpha|}
-\sum_{\substack{t,s\in\mathbb N_0\\t+s=|\alpha|-1}}
\left(1+\frac{x_0}{\vec{q}}\right)^t\frac{x_0}{\vec{q}}
\left(1+\frac{x_0}{\vec{q}}\right)^s\right]\\
&\hspace{-4cm}=
\frac{x^\alpha}{\vec{q}}\left[|\alpha|\left(1+\frac{x_0}{\vec{q}}\right)^{|\alpha|}
-|\alpha|\left(1+\frac{x_0}{\vec{q}}\right)^{|\alpha|-1}\frac{x_0}{\vec{q}}\right].
    \end{split}
  \]

  STEP 5: {\rm We can now compute $V_q\mu^\alpha$.}\\

 Using \eqref{osaka12345} and the previous step we have:
  \[
\begin{split}
    V_q\mu^\alpha&=    x^\alpha\frac{|\alpha|}{\vec{q}}  \left(1+\frac{x_0}{\vec{q}}\right)^{|\alpha|-1}-
\frac{x^\alpha}{\vec{q}}\left[|\alpha|\left(1+\frac{x_0}{\vec{q}}\right)^{|\alpha|}
  -|\alpha|\left(1+\frac{x_0}{\vec{q}}\right)^{|\alpha|-1}\frac{x_0}{\vec{q}}\right]\\
&=\frac{|\alpha|x^\alpha}{\vec{q}}\left(1+\frac{x_0}{\vec{q}}\right)^{|\alpha|-1}
\left(  1-\left(1+\frac{x_0}{\vec{q}}\right)+\frac{x_0}{\vec{q}}\right)\\
&=0.
\end{split}
\]
This ends the proof.
 \end{proof}

 \begin{remark}
We observe that both functions $\mu^\alpha$ and $\zeta^\alpha$ coincide with $x^\alpha$ when $x_0=0$. It is important to note that these are two different extensions of the same real function $x^\alpha$ leading to two different regular function theories. In fact, $\mu^\alpha$ is the $V_q$-regular extension of $x^\alpha$ while $\zeta^\alpha$ gives the classical Fueter extension. However, the classical Fueter variables $\zeta^\alpha$ extend $x^\alpha$ to the whole space of quaternions while $\mu^\alpha$ extend $x^\alpha$ to domains of $\mathbb{H}^*$.     
    \end{remark}

  \begin{proposition}
   We prove that $q^n$ is in $\ker V_q$ and that, moreover ,
    \[
q^n=\sum_{|\alpha|=n}\mu^\alpha c_{\alpha,n}
\]
where, with $\alpha=(\alpha_1,\alpha_2,\alpha_3)$,
\begin{equation}
c_{\alpha,n}=\frac{n!}{\alpha!}{\mathbf e_1}^{\times \alpha_1}\times{\mathbf e_2}^{\times \alpha_2}\times{\mathbf e_3}^{\times \alpha_3},
\label{cocosisisicocosi}
\end{equation}
where the symmetric product is taken among all the products of the units $\mathbf{e_u}$.
  \end{proposition}

  \begin{proof}
  In $\mathbb{H}^*$ we have 
    \[
  \begin{split}
    q^n&=(x_0+\vec{q})^n\\
    &=\left(1+\frac{x_0}{\vec{q}}\right)^n  \left(\vec{q}\,\right)^n\\
    &=\left(1+\frac{x_0}{\vec{q}}\right)^n  \left(\sum_{|\alpha|=n} x^\alpha c_{\alpha,n}\right)\\
      &=\sum_{|\alpha|=n}x^\alpha\left(1+\frac{x_0}{\vec{q}}\right)^n c_{\alpha,n},
    \end{split}
  \]
  for some $c_{\alpha,n}\in\mathbb H$ which can be expressed in term of symmetrized products as in
  \eqref{cocosisisicocosi}
    by known formulas.
\end{proof}

We note that \eqref{cocosisisicocosi} does not take into account  the Cayley table of multiplication for the quaternions.




\begin{example}
  \label{totoche}
Let us examine $\displaystyle\mu_1(x)=x_1\left(1+\frac{x_0}{\vec{q}}\right)=x_1+x_0\frac{x_1}{\vec{q}}$. The function
 $\displaystyle\frac{\partial\mu_1}{\partial x_2}=x_1x_0\left(\frac{\mathbf e_2}{\vec{q}^{\, 2}}+\frac{2x_2}{\vec{q}^{\, 3}}\right)$ is not $V_q$-regular, while the function
 $\displaystyle\frac{\partial\mu_1}{\partial x_0}=\frac{x_1}{\vec{q}}$ is $V_q$ -regular but cannot
 be written as a convergent Fueter-like series of the $\mu^\alpha$.

\end{example}

\begin{proof}[Discussion of Example \ref{totoche}]
  We have
  \[
\frac{\partial \mu_1}{\partial x_2}(x)=x_1x_0\left(\frac{\mathbf e_2}{\vec{q}^{\, 2}}+\frac{2x_2}{\vec{q}^{\, 3}}\right),
    \]
and so
\[
\frac{\partial^2 \mu_1}{\partial x_0 \partial x_2}(x)=x_1\left(\frac{\mathbf e_2}{\vec{q}^{\, 2}}+\frac{2x_2}{\vec{q}^{\, 3}}\right)
\]
is independent of $x_0$. On the other hand
\[
  \begin{split}
    \sum_{u=1}^3x_u\frac{\partial}{\partial x_u}\left(\frac{\partial \mu_1}{\partial x_2}\right)&=
    \sum_{u=1}^3x_u
  \frac{\partial}{\partial x_u}\left(x_1x_0\left(\frac{\mathbf e_2}{\vec{q}^{\, 2}}+\frac{2x_2}{\vec{q}^{\, 3}}\right)\right)\\
  &=x_0\left[\sum_{u=1}^3x_u
  \frac{\partial}{\partial x_u}\left(x_1\left(\frac{\mathbf e_2}{\vec{q}^{\, 2}}+\frac{2x_2}{\vec{q}^{\, 3}}\right)\right)\right].
\end{split}
\]
Hence
\[
\left(V_q\frac{\partial \mu_1}{\partial x_2}\right)(x)=x_1\left(\frac{\mathbf e_2}{\vec{q}^{\, 2}}+\frac{2x_2}{\vec{q}^{\, 3}}\right)-
x_0\left[\sum_{u=1}^3x_u
  \frac{\partial}{\partial x_u}\left(x_1\left(\frac{\mathbf e_2}{\vec{q}^{\, 2}}+\frac{2x_2}{\vec{q}^{\, 3}}\right)\right)\right].
\]
Setting $x_0=0$ shows that $\displaystyle V_q\frac{\partial \mu_1}{\partial x_2}\not\equiv 0$.\\

We now turn to  $\displaystyle\frac{\partial\mu_1}{\partial x_0}$, which  is $V_q$ -regular since differentiation with respect to $x_0$
commutes with $V_q$.  Assume now that
  \[
\frac{\partial\mu_1}{\partial x_0}=\sum_{\alpha\in\mathbb N_0^3}\mu^\alpha c_\alpha
\]
for some quaternionic coefficients $c_\alpha$, and where the convergence is assumed in some
set $\Omega\subset \mathbb{H}^*$. In other words
\begin{equation}
\frac{x_1}{\vec{q}}=\sum_{\alpha\in\mathbb N_0^3}x^\alpha\left(1+\frac{x_0}{\vec{q}}\right)^{|\alpha|}c_\alpha.
\label{atlantique}
\end{equation}
Setting $x_0=0$ in this equation, and multiplying both sides by $\vec{q}$ leads to
\begin{equation}
\label{atlantique2}
x_1=\sum_{\alpha\in\mathbb N_0^3}x_1^{\alpha_1+1}x _2^{\alpha_2}x_3^{\alpha_3}{\mathbf e_1}c_\alpha+
\sum_{\alpha\in\mathbb N_0^3}x_1^{\alpha_1}x _2^{\alpha_2+1}x_3^{\alpha_3}{\mathbf e_2}c_\alpha+
\sum_{\alpha\in\mathbb N_0^3}x_1^{\alpha_1}x _2^{\alpha_2}x_3^{\alpha_3+1}{\mathbf e_3}c_\alpha,
\end{equation}
where this equality is valid {\sl a priori} only for $(x_1,x_2,x_3)$ such that $x_u\neq 0$ for some $u=1,2,3$. By classical results on
series and summable families \eqref{atlantique2} can be extended to $x_u=0$.  Comparing the linear terms
on both sides of \eqref{atlantique2}  leads to
\[
  x_1=x_1{\mathbf e_1}c_{0,0,0}+x_2{\mathbf e_2}c_{0,0,0}+x_3{\mathbf e_3}c_{0,0,0}
\]
which is impossible, leading to a contradiction.
\end{proof}

Just as in the classical case, let us now define the conjugate operator of $V_q$, denoted by $\overline{V_q}$, to be the operator defined by: 
\begin{equation}
\displaystyle \overline{V_q}:=\frac{\partial}{\partial x_0}+\frac{1}{\vec{q}}\sum_{u=1}^{3}\mathbf{e}_u\frac{\partial}{\partial x_u}.
\end{equation}
We can prove the following:
\begin{proposition}
Let $\Omega$ be an open domain in $\mathbb{H}^*$.
If $f$ is $V_q$-regular on $\Omega$, then $\overline{V_q}(f)$ is also $V_q$-regular on $\Omega$. Moreover, in this case we have 
\begin{equation}
\displaystyle \frac{1}{2}\overline{V_q}f(q)=\frac{\partial}{\partial x_0}f(q), \quad \forall q=x_0+\vec{q}.
\end{equation}
\end{proposition}
\begin{proof}
We first observe that $$\displaystyle (V_q+\overline{V_q})f(q)=2\frac{\partial }{\partial x_0}f(q).$$

However, since $f$ is $V_q$-regular, then it belongs to $\ker(V_q)$ so that $V_q(f)=0$, leading to 
$$\displaystyle \overline{V_q}f(q)=2\frac{\partial}{\partial x_0}f(q).$$
Thus, applying the operator $V_q$ and, using the fact that it commutes with $\displaystyle\frac{\partial}{\partial x_0}$, we obtain \begin{align*}
\displaystyle V_q\overline{V_q}f(q)&=2V_q\frac{\partial}{\partial x_0 }f(q)\\
&=2\frac{\partial}{\partial x_0}V_qf(q)\\
&=0.\\
\end{align*}
So, the function $\overline{V_q}(f)$ is $V_q$-regular, which ends the proof.
\end{proof}

\begin{proposition}
It holds that for every $q\in\mathbb{H}^*$ we have
$$\displaystyle \frac{1}{2} \overline{V_q}\mu^\alpha(q)=\frac{|\alpha|}{q}\mu^\alpha(q).$$
\end{proposition}
\begin{proof}
We use the previous computations of $\displaystyle\frac{\partial}{\partial x_0}\mu^\alpha$ and $\displaystyle\frac{1}{\vec{q}}\sum_{u=1}^{3}x_u\frac{\partial}{\partial x_u}\mu^\alpha$ to obtain 
$$\displaystyle \overline{V_q}\mu^\alpha(q)=2|\alpha|\frac{x^\alpha}{\vec{q}}\left(1+\frac{x_0}{\vec{q}}\right)^{|\alpha|-1}.$$
Hence, we get 

$$\displaystyle \overline{V_q}\mu^\alpha(q)=2|\alpha|q^{-1}\mu^\alpha(x),$$
which is the desired equality.
\end{proof}

\begin{remark}
It is easy to see that a consequence of this proposition yields the counterpart of the properties of the number operators and we have 
$$M_q\overline{V_q}(\mu^{\alpha})=2|\alpha| \mu^{\alpha},$$ 
where $M_q$ represents the left quaternionic multiplication.
\end{remark}

\section{Cauchy-Kovalevskaia product}
\setcounter{equation}{0}
\label{CK-Vq}

In this section we build a Cauchy-Kovalevskaia product using the new variables $\mu^\alpha$.
We start by building a CK-extension of a function in the kernel of $V_q$, $f=f_0+\mathbf e_1f_1+\mathbf e_2f_2+\mathbf e_3 f_3$, where $f_0,f_1,f_2,f_3$ are real valued differentiable functions on an open subset $\Omega \subset \mathbb{H}^*$.
The equation $V_qf=0$ can be rewritten as
\begin{equation}
  \label{expre}
\frac{\partial  f_j}{\partial x_0}=G_j,\quad j=0,1,2,3\, ,
\end{equation}
where $G_j$ is real analytic in the variables $x_0,x_1,x_2,x_3$. The right hand side can be written in terms of the partial derivatives with respect to the other three variables, i.e. $\displaystyle\frac{\partial f_j}{\partial x_k},\,j=0,1,2,3,\, k=1,2,3$,
and they are given by the following system on an open subset $\Omega\subset\mathbb{H}^*$
\begin{equation}
\begin{split}
\frac{\partial f_0}{\partial x_0}&=\frac{1}{x_1^2+x_2^2+x_3^2}\left(
  x_1E f_1+x_2Ef_2+x_3Ef_3\right)\\
\frac{\partial f_1}{\partial x_0}&=-\frac{1}{x_1^2+x_2^2+x_3^2}\left(-x_1Ef_0-x_2Ef_3+x_3Ef_2\right)\\
\frac{\partial f_2}{\partial x_0}&=-\frac{1}{x_1^2+x_2^2+x_3^2}\left(-x_2Ef_0-x_1Ef_3+x_3Ef_1\right)\\
\frac{\partial f_3}{\partial x_0}&=-\frac{1}{x_1^2+x_2^2+x_3^2}\left(-x_3Ef_0-x_1Ef_2+x_2Ef_1\right)
\end{split}
\end{equation}
where $E$ denotes the Euler operator:
\[
E=\sum_{u=1}^3x_u\frac{\partial}{\partial x_u}.
  \]
We can therefore apply the
Cauchy-Kovalevskaia theorem (see e.g. \cite[\S7 p. 39]{MR25:4216}, which asserts that the system \eqref{expre} has a
unique solution near a real point possibly different from the origin) for given initial real analytic values
$f_j(0,x_1,x_2,x_3)$ on an open domain $\widetilde{\Omega}\subset (\mathbb{H}^*\cap\,\{x_0=0\})$. This solution, $F$, is defined on an open set $\Omega\subset\mathbb{H}^*$, where $\widetilde{\Omega}=\Omega\, \cap \,\{x_0=0\}$ and we call it the {\em CK-extension of $f$ with respect to $V_q$}. 

\begin{definition} The CK-extension to $\Omega$ with respect to the operator $V_q$ found above is denoted by $F=CK_{V_q}(f)$.
  \end{definition}

\begin{remark} The choice
\[
f_0(0,x_1,x_2,x_3)=x_1^{\alpha_1}x_2^{\alpha_2}x_3^{\alpha_3}\quad{\rm and}\quad f_1(0,x_1,x_2,x_3)=f_2(0,x_1,x_2,x_3)=f_3(0,x_1,x_2,x_3)=0
\]
leads, in the case of the variables $x^\alpha$ to:
  \begin{equation}
    \mu^\alpha=CK_{V_q}(x^\alpha).
    \end{equation}
\end{remark}

\begin{definition} We can now define the CK-product of two functions $f, g$ in the kernel of $V_q$ to be:
  \begin{equation}
    f\star_{V_q} g=CK_{V_q}(f(0,x_1,x_2,x_3)g(0,x_1,x_2,x_3)),
    \end{equation}
 for every $x\in\Omega$.
\end{definition}

It is easy to see that, for any two quaternions $c,d$ we have:
\begin{proposition}  On the entire domain $\mathbb{H}^*$ we have that 
  \begin{equation}
\mu^\alpha c\star_{V_q}\mu^\beta d=\mu^{\alpha+\beta}cd.
  \end{equation}
\end{proposition}

\begin{proof}
The proof is left to the reader.
  \end{proof}

  \begin{theorem}   For every $x\in\mathbb{H}^*$ it holds that:
    \begin{equation}
      \label{okinawa}
\left(\exp\left(x_0\frac{1}{\vec{q}}\,E\right)\right)(x^\alpha)=\mu^\alpha.
      \end{equation}
  \end{theorem}

  \begin{proof} We proceed in a number of steps.\\

    STEP 1: {\sl It holds that
      \begin{equation}
      \label{by}
        E\left(\frac{1}{\vec{q}}\right)=-\frac{1}{\vec{q}} \, .
      \end{equation}
    }

    Indeed, using \eqref{kobe} we can write
\[
    \begin{split}
      E\left(\frac{1}{\vec{q}}\right)&=\sum_{u=1}^3x_u\left(\frac{\mathbf e_u}{\vec{q}^{\, 2}}+\frac{2x_u}{\vec{q}^{\, 3}}\right)\\
      &=\frac{\vec{q}}{\vec{q}^{\,2}}+2\frac{x_1^2+x_2^2+x_3^2}{\vec{q}^{\, 3}}\\
  &=\frac{\vec{q}}{\vec{q}^{\, 2}}-2\frac{\vec{q}^{\,2}}{\vec{q}^{\, 3}}\\
  &=-\frac{1}{\vec{q}}\,.
    \end{split}
  \]
  
    STEP 2: {For non-commuting functions of a real variable it holds that}
    \begin{equation}
      \label{vienne}
      (fg)^\prime=f^\prime g+fg^\prime.
      \end{equation}
      This follows from
      \[
        \begin{split}
f(t)g(t)-f(t_0)g(t_0)&=(f(t)-f(t_0))g(t)+f(t_0)(g(t)-g(t_0)).
        \end{split}
        \]
Note that $f$ and $g$ may commute at a joint value $t$ but we do not assume that $f(t)g(s)=g(s)f(t)$ for $t\not=s$.\\
      
    STEP 3: {\sl It holds that
      \begin{equation}E\left(\frac{x^\alpha}{\vec{q}^{\, n}}\right)=\begin{cases} (|\alpha|-n)\frac{x^\alpha}{\vec{q}^{\, n}},\quad 0\le n\le
          |\alpha|,\\
            \, 0,\quad\hspace{1.75cm} otherwise.
          \end{cases}
          \label{recur567}
        \end{equation}
      }
      We proceed by induction. The case $n=0$ corresponds to the formula
      \[
        E(x^\alpha)=|\alpha|x^\alpha.
        \]
        We then write
      \[
        \frac{x^\alpha}{\vec{q}^{\, n+1}}=        \frac{x^\alpha}{\vec{q}^{\, n}}\frac{1}{\vec{q}}
        \]
      and apply \eqref{vienne} with $\displaystyle f(x)=\frac{x^\alpha}{\vec{q}^{\, n}}$ and $\displaystyle g(x)=\frac{1}{\vec{q}}$ to get:
      \[
        \begin{split}
          E\left(\frac{x^\alpha}{\vec{q}^{\,n+1}}\right)&=\sum_{u=1}^3x_u\frac{\partial}{\partial x_u}\left(\frac{x^\alpha}{\vec{q}^{n+1}}\right)\\
         &=\left(\sum_{u=1}^3   x_u\frac{\partial}{\partial x_u}\left(     \frac{x^\alpha}{\vec{q}^{\, n}}\right)\right)\frac{1}{\vec{q}}+
          \frac{x^\alpha}{\vec{q}^{\, n}}\left(\sum_{u=1}^3x_u\frac{\partial}{\partial x_u}\frac{1}{\vec{q}}\right)\\
          &=\underbrace{(|\alpha|-n)\frac{x^\alpha}{\vec{q}^{\, n}}}_{\mbox{{\rm induction at rank $n$}}}
          \cdot \frac{1}{\vec{q}}+\frac{x^\alpha}{\vec{q}^{\, n}}\underbrace{\frac{-1}{\vec{q}}}_{\mbox{\rm by~(\ref{by})}}\\
          &=(|\alpha|-n-1)\frac{x^\alpha}{\vec{q}^{\, n+1}}.
            \end{split}
          \]

          STEP 4: {\sl We prove
            \begin{equation}
              \label{toulon}
              \left(E\, \frac{1}{\vec{q}}\right)^n(x^\alpha)=(|\alpha|-1)(|\alpha|-2)\cdots (|\alpha|-n+1)\frac{x^\alpha}{\vec{q}^{\, n}},
              \quad n=1,2,\ldots
              \end{equation}
            }
            We proceed by induction. The case $n=1$ corresponds to the previous step.
            \[
              \begin{split}
          \left(E\, \frac{1}{\vec{q}}\right)^{n+1}(x^\alpha)&=E\,\frac{1}{\vec{q}}\left(\left(E\,\frac{1}{\vec{q}}\right)^n(x^\alpha)\right)
                  \\
                  &=E\,\frac{1}{\vec{q}}\left(
                    (|\alpha|-1)(|\alpha|-2)\cdots (|\alpha|-n+1)\frac{x^\alpha}{\vec{q}^{\, n}}                  \right)\\
                  &=(|\alpha|-1)(|\alpha|-2)\cdots (|\alpha|-n+1)
                  E\left(\frac{x^\alpha}{\vec{q}^{\, n}}\right)\frac{1}{\vec{q}}\\
                  &=(|\alpha|-1)(|\alpha|-2)\cdots (|\alpha|-n+1)(|\alpha|-n) \frac{x^\alpha}{\vec{q}^{\, n+1}},
                \end{split}
              \]
 where we have used the induction hypothesis to go from the first to the second line, and  formula \eqref{recur567} to go from
 the third line to the fourth.\\
            
 STEP 5:{\sl\, We prove \eqref{okinawa}.}\smallskip

 We first note that, in view of \eqref{toulon}
 \[
   \left(   E\, \frac{1}{\vec{q}}\right)^n(x^\alpha)=0
 \]
 for $n\ge |\alpha |$. Furthermore, since
 \[
\left( \frac{1}{\vec{q}} \, E\right)^n= \frac{1}{\vec{q}} \left(   E\, \frac{1}{\vec{q}}\right)^{n-1}E,\quad n=1,2,\ldots,
\]
we have that the series
 \[
   \begin{split}
     \left(\exp\left(x_0\, \frac{1}{\vec{q}}\, E\right)\right)(x^\alpha)&=x^\alpha+\sum_{n=1}^\infty\frac{x_0^n}{n!}
     \left(\frac{1}{\vec{q}}\, E\right)^n (x^\alpha)\\
     &=x^\alpha+\sum_{n=1}^\infty\frac{x_0^n}{n!}\frac{1}{\vec{q}}
     \left(E\, \frac{1}{\vec{q}}\right)^{n-1}E (x^\alpha)\\
   \end{split}
   \]
has only a finite number of non-zero terms, and is therefore equal to
 \[
   \begin{split}
     \left(\exp\left(x_0\, \frac{1}{\vec{q}}\, E\right)\right)(x^\alpha)
     &=x^\alpha+\sum_{n=1}^\infty\frac{x_0^n}{n!}\frac{1}{\vec{q}}
     \left(E\, \frac{1}{\vec{q}}\right)^{n-1}E (x^\alpha)\\
     &=x^\alpha\left(\sum_{n=0}^{|\alpha|-1}\frac{|\alpha|(|\alpha|-1)\cdots(|\alpha|-n+1)}{n!}
       \left(\frac{x_0}{\vec{q}}\right)^n\right)\\
     &=x^\alpha\left(1+\frac{x_0}{\vec{q}}\right)^{|\alpha|}.
   \end{split}
   \]
The theorem is now proven.

  \end{proof}

  Let $b$ and $c$ be two quaternions and writing
    \[
    q^nb=\sum_{\alpha\in\mathbb N_0^3}\mu^\alpha c_{\alpha,n}b\quad{\rm and}\quad
  q^mc=\sum_{\alpha\in\mathbb N_0^3}\mu^\alpha c_{\alpha,m}c,
\]
we have the $\star_{V_q}$ product of the two:
\[
  q^nb\star_{V_q}q^mc=\sum_{\alpha,\beta\in\mathbb N_0^3}\mu^{\alpha+\beta}c_{\alpha,n}bc_{\beta,m}c.
\]
In the expression above, since $b$ and $c$ are quaternions they do not commute with $c_{\alpha,n}$ and $c_{\beta,m}$, however, when $b\in\mathbb R$ this expression reduces to
\begin{equation}
  \begin{split}
    q^nb\star_{V_q}q^mc&=\sum_{\alpha,\beta\in\mathbb N_0^3}\mu^{\alpha+\beta}c_{\alpha,n}c_{\beta,m}bc\\
    &=\sum_{\gamma\in\mathbb N_0^3}\left(\sum_{\substack{\alpha+\beta=\gamma\\ \alpha,\beta\in\mathbb N_0^3}}c_{\alpha,n}c_{\beta,m}\right)bc,\quad b\in\mathbb R,\,\,
    c\in\mathbb H.
    \end{split}
  \end{equation}

  When $x_0=0$ we have $q=\vec{q}$ and the $\star_{V_q}$ product reduces to the pointwise product.
  We have
  \[
\vec{q}^{\, \, n}=\sum_{\alpha\in\mathbb N_0}x^\alpha c_{\alpha,n}
    \]
    and
    \[
      \vec{q}^{\,\, n}\vec{q}^{\,\, m}=\vec{q}^{\, \, n+m}
      \]
    so that
    \[
c_{\gamma,n+m}=\sum_{\substack{\alpha+\beta=\gamma\\ \alpha,\beta\in\mathbb N_0^3}}c_{\alpha,n}c_{\beta,m}.
      \]
  \begin{equation}
  \begin{split}
  q^nb\star_{{\rm slice}}q^mc&=q^{n+m}bc\\
  &=\sum_{\gamma\in\mathbb N_0^3}\mu^{\gamma}c_{\gamma,n+m}bc.
    \end{split}
  \end{equation}

We then have:
  
  \begin{proposition}

    Let $b\in\mathbb R$ and $c\in\mathbb H$. On the entire domain $\mathbb{H}^*$ it holds that
    \begin{equation}
q^nb\star_{V_q}q^mc=q^nb\star_{{\rm slice}}q^mc=q^{n+m}bc
\end{equation}
and, more generally, for $f$ intrinsic
\[
f\star_{V_q}g=f\star_{{\rm slice}}g=fg.
  \]
    \end{proposition}
  
  The term {\em intrinsic} used in slice quaternionic analysis means that the function preserves all slices.
  This equality does not hold for more general functions in the kernel of $V_q$.

  \section{Reproducing kernel Hilbert spaces}
  \label{Vq-Kernel-Omega}
  \setcounter{equation}{0}

We now set the stage for the Schur analysis theory in this case and we begin by defining a reproducing kernel for the Hilbert space of power series in $\mu^\alpha$.\\
For $r,R$ and $\rho$ strictly positive, let us first define the domain:
\begin{equation}
 \Omega_{r,R, \rho}=\left\{x\in\mathbb R^4\,:\,  r<|x_u|<R,\,\, u=1,2,3,\,\,\,{\rm and},\,\, |x_0|<\rho\right\}.
  \end{equation}
  
 First we see that for any $x\in\Omega_{r,R,\rho}\subset \mathbb{H}^*$, we have:
\begin{equation}
  |\mu^\alpha(x)|\le   L_{r,R,\rho}^{|\alpha|},
\end{equation}
where $\displaystyle L_{r,R,\rho}=R\left(1+\frac{\rho}{\sqrt{3}r}\right)$.

Let $\Omega_{r,R,\rho}$ and $L$ be defined as above.

\begin{proposition} Let $c_{\alpha}$ be a family of positive numbers for $\alpha$ in a subset $S$ of $\mathbb N_0^3$. Assume that
  \[
\sum_{\alpha\in S} \frac{L_{r_1,R_1,\rho_1}^{2|\alpha|}}{c_\alpha}<\infty
\]
for all $r_1,R_1,\rho_1$ such that
\begin{equation}
  \label{toyota}
  r<r_1<R_1<R  \quad and\quad 0<\rho_0<\rho.
  \end{equation}
Then the function
  \begin{equation}
K_{\bf c}(x,y)=\sum_{\alpha\in S}\frac{\mu^\alpha(x)\overline{\mu^\alpha(y)}}{c_\alpha},
\end{equation}
is positive definite in $\Omega_{r,R,\rho}$ and the associated reproducing kernel Hilbert space consists of the Fueter-like series with quaternionic coefficients $f_{\alpha}$:

\begin{equation}
\mathfrak H(K_{\mathbf c})=\{f= \sum_{\alpha\in\mathbb N_0^3}\mu^{\alpha}\,f_\alpha\,\, | \,\, f \, {\rm abs. \,  conv.\,  in\,}\Omega_{r,R,\rho},\, \sum_{\alpha\in\mathbb N_0^3}c_\alpha |f_\alpha|^2<\infty\}.
\end{equation}

  \end{proposition}

  \begin{proof}
The proof follows from the usual arguments as in~\cite{MR3819695, MR3795166}.
  \end{proof}

  
  \begin{proposition}
Elements of $\mathfrak H(K_{\mathbf c})$ are $V_q$-regular on
    \begin{equation}
\Omega(\mathfrak H(K_{\mathbf c}))=\left\{x\in\mathbb{H}^*\,|\, \sum_{\alpha\in\mathbb N_0^3}\frac{|\mu^\alpha(x)|^2}{c_\alpha}<\infty \right\}.
      \end{equation}
   
    \end{proposition}

     \begin{proof} We proceed in a number of steps.\\

       STEP 1: {\sl Let  $a\in(0,1)$. Then
         \[
           \sum_{\alpha\in\mathbb N_0^3}a^{2|\alpha|}|\alpha|^2<\infty,
         \]
         and in particular $a^{|\alpha|}|\alpha|$ are uniformly bounded.
       }
       By the Cauchy-Schwarz inequality
       \[
         ( \alpha_1+\alpha_2+\alpha_3)^2\le 3(\alpha_1^2+\alpha_2^2+\alpha_3^2).
       \]
       Hence
       \[
         \begin{split}
           \sum_{\alpha\in\mathbb N_0^3}a^{2|\alpha|}|\alpha|^2&\leq
           3\sum_{\alpha\in\mathbb N_0^3}a^{2|\alpha|}(\alpha_1^2+\alpha_2^2+\alpha_3^2)\\
           &=3\sum_{u=1}^3\left(\sum_{\alpha\in\mathbb N_0^3}a^{2\alpha_1}a^{2\alpha_2}a^{2\alpha_3}\alpha^2_u\right).
\end{split}
        \] 
        But
        \[
          \left(\sum_{\alpha\in\mathbb N_0^3}a^{2\alpha_1}a^{2\alpha_2}a^{2\alpha_3}\alpha_1^2\right)=
            \left(\sum_{\alpha_1=0}^\infty a^{2\alpha_1}\alpha_1^2              \right)
            \left(\sum_{\alpha_2=0}^\infty a^{2\alpha_2}              \right)
            \left(\sum_{\alpha_3=0}^\infty a^{2\alpha_3}             \right)<\infty
          \]
          and similarly for $u=2$ and $u=3$.\\
          
       STEP 2: {\sl For $r_1$ and $R_1$ as in \eqref{toyota}
         \begin{equation}
           \label{wertyu2}
\sum_{\alpha\in\mathbb N_0^3}\frac{\left|\frac{\partial \mu^\alpha}{\partial x_1}\right|^2}{c_\alpha}<\infty
\end{equation}
for $r_1\le |x_1| \le R_1$,} and similarly for $x_2$ and $x_3$.\smallskip

We have
\[
  \begin{split}
    \frac{\partial \mu^\alpha}{\partial x_1}&=\alpha_1x^{\alpha-e_1}\left(1+\frac{x_0}{\vec{q}}\right)^{|\alpha|}+\\
    &\hspace{5mm}+x^\alpha\left[\sum_{\substack{t,s\in\mathbb N_0\\ t+s=|\alpha|-1}}
        \left(1+\frac{x_0}{\vec{q}}\right)^t\left(\frac{x_0\mathbf e_u}{\vec{q}^{\,2}}+\frac{2x_ux_0}{\vec{q}^{\,3}}\right)
        \left(1+\frac{x_0}{\vec{q}}\right)^s\right]\\
\end{split}
\]
and so
\[
  \begin{split}
\left|\frac{\partial \mu^\alpha}{\partial x_1}\right|&\le\alpha_1R_1^{|\alpha|-1}+R_1^{|\alpha|}|\alpha|
\left(1+\frac{\rho_1}{\sqrt{3}r_1}\right)^{|\alpha|-1}
\left(\frac{\rho_1}{3r_1^2}+\frac{2R_1\rho_1}{3^{3/2}r_1^{3/2}}\right)\\
&\le R_1^\alpha\left(1+\frac{\rho_1}{\sqrt{3}r_1}\right)^{|\alpha|}|\alpha|
\underbrace{\left[1+\frac{1}{1+\frac{\rho_1}{\sqrt{3}r_1}}
      \left(\frac{\rho_1}{3r_1^2}+\frac{2R_1\rho_1}{3^{3/2}r_1^{3/2}}\right)
 \right]}_{\mbox{=$M$, \rm independent of $\alpha$}}.
  \end{split}
\]
Let $R_2$ be such that $R_1<R_2<R$. We can write
\[
  \begin{split}
    \sum_{\alpha\in\mathbb N_0^3}\frac{\left|\frac{\partial \mu^\alpha}{\partial x_1}\right|^2}{c_\alpha}&\le \, M
    \sum_{\alpha\in\mathbb N_0^3}\frac{   R_1^{2|\alpha|}\left(1+\frac{\rho_1}{\sqrt{3}r_1}\right)^{2|\alpha|}|\alpha|^2 }{c_\alpha} \\
    &= \, M\sum_{\alpha\in\mathbb N_0^3}\frac{   R_2^{2|\alpha|}\left(1+\frac{\rho_1}{\sqrt{3}r_1}\right)^{2|\alpha|}|\alpha| }{c_\alpha}
    \left[\left(
        \frac{R_1}{R_2}\right)^{2|\alpha|}|\alpha|^2\right] \\
    &\le \, M \sum_{\alpha\in\mathbb N_0^3}\frac{   R_2^{2|\alpha|}\left(1+\frac{\rho_1}{\sqrt{3}r_1}\right)^{2|\alpha|}|\alpha| ^2}{c_\alpha}\\
    &<\infty\, 
  \end{split}
\]
using Step 1, for some constant $M$.\\

STEP 3: {\sl For $r_1$ and $R_1$ as in \eqref{toyota} we have:
         \begin{equation}
           \label{wertyu3}
\sum_{\alpha\in\mathbb N_0^3}\frac{\left|\frac{\partial \mu^\alpha}{\partial x_0}\right|^2}{c_\alpha}<\infty
\end{equation}
for $|x_0|\le \rho_1$.}\\ 

This follows from
\[
  \frac{\partial}{\partial x_0}\mu^\alpha=x^\alpha
  \left[\sum_{\substack{t,s\in\mathbb N_0\\ t+s=|\alpha|-1}}
        \left(1+\frac{x_0}{\vec{q}}\right)^t\left(\frac{1}{\vec{q}}\right)
        \left(1+\frac{x_0}{\vec{q}}\right)^s\right]
    \]
   and the corresponding bound 

\[
  \left|\frac{\partial}{\partial x_0}\mu^\alpha\right|\le |\alpha|R_1^{|\alpha|}\left(1+\frac{\rho_1}{\sqrt{3}r_1}\right)^{|\alpha|-1}\frac{1}{\sqrt{3}r_1}.
\]
\mbox{}\\

This step can be proven directly using the Appell-type property, i.e. $V_{\overline{q}}$ acts on powers of $\mu$ the same as $\displaystyle \frac{\partial}{\partial x_0}$.\\

STEP 4: {\sl We prove that, pointwise, for $f=\sum_{\alpha\in\mathbb N_0^3}\mu^\alpha f_\alpha\in\mathfrak H(K_{\mathbf c})$,
  \[
    \frac{\partial}{\partial x_u}\sum_{\alpha\in\mathbb N_0^3}\mu^\alpha f_\alpha
        =
\sum_{\alpha\in\mathbb N_0^3}\frac{\partial}{\partial x_u}\mu^\alpha f_\alpha .
\]
}
Using the Cauchy-Schwarz inequality and the previous lemma we see that the series of derivatives
\[
    \sum_{\alpha\in\mathbb N_0^3}\frac{\partial}{\partial x_u}
  \mu^\alpha f_\alpha=
  \sum_{\alpha\in\mathbb N_0^3}\frac{\frac{\partial}{\partial x_u}\mu^\alpha }{\sqrt{c_\alpha}}\sqrt{c_\alpha}f_\alpha
 \]
 converges uniformly in intervals $[r_1,R_1]$.
 Since the series $\sum_{\alpha\in\mathbb N_0^3}\mu^\alpha f_\alpha$ converges in $\mathfrak H(K_{\mathbf c})$ then it converges pointwise and a classical calculus theorem allows us to conclude that we have convergence in Step 4. 
 This classical theorem speaks of sequences and not of summable families, but we can reduce the latter to the case of sequences by identifying $\mathbb N_0^3$ and $\mathbb N_0$ via a bijection.\\

 This solves the cases of $u=1,2,3$. The case $u=0$ is treated in a similar way and we leave the details to the reader.\\



We now apply Step 4 four times in the definition of $V_q$. More precisely,
 \[
   \begin{split}
   V_q\left(\sum_{\alpha\in\mathbb N_0^3}\mu_\alpha f_\alpha\right)
   &=\frac{\partial}{\partial x_0}\left(\sum_{\alpha\in\mathbb N_0^3}\mu_\alpha f_\alpha\right)
   -\frac{1}{\vec{q}}\sum_{u=1}^3
   x_u\frac{\partial}{\partial x_u}\left(\sum_{\alpha\in\mathbb N_0^3}\mu_\alpha f_\alpha\right)\\
   &=   \left(\sum_{\alpha\in\mathbb N_0^3}\frac{\partial}{\partial x_0}\mu_\alpha f_\alpha\right)-\\
   &\hspace{5mm}-\frac{1}{\vec{q}}\sum_{u=1}^3
   x_u   \left(\sum_{\alpha\in\mathbb N_0^3}\frac{\partial}{\partial x_u}\mu_\alpha f_\alpha\right)\\
   &=\sum_{\alpha\in\mathbb N_0^3} V_q\mu^\alpha f_\alpha\\
   &=0.
     \end{split}
   \]
This concludes the proof.

  \end{proof}
Let us now turn to a definition of a backward-shift and multiplication operator in this case and we set
\begin{equation}
  \label{leipzig123}
  e_1=(1,0,0), \quad e_2=(0,1,0)\quad{\rm  and}\quad e_3=(0,0,1).
  \end{equation}

  \begin{definition}

   The multiplication and backward shift operators on $\mathbb{H}^*$ are:
    \begin{eqnarray}
      M_u\mu^\alpha&=&\mu^{\alpha+e_u}\\
      B_u\mu^\alpha&=&\frac{\alpha_u}{|\alpha|}\mu^{\alpha-e_u},\quad u=1,2,3
    \end{eqnarray}
    with the understanding that $B_u\mu^\alpha=0$ if $\alpha_u=0$.
   
  \end{definition}

  Note that
  \begin{equation}
\left(\sum_{u=1}^3B_uM_u\right)\mu^\alpha=\mu^\alpha.
\end{equation}

As expected, we have the following result:
  
  \begin{theorem}
    Assume that the domain of definition of $M_u$ in $\mathfrak H(K_{\mathbf c})$ contains the linear span of the Fueter polynomials $\mu^\alpha$.
    Then, 
    \begin{equation}
      \label{alpha-u}
M_u^*(\mu^\alpha)=\frac{c_\alpha}{c_{\alpha-e_u}}\mu^{\alpha-e_u},
\end{equation}
with the understanding that the right hand side of \eqref{alpha-u} equals $0$ if $\alpha_u=0$.
    \end{theorem}

    \begin{proof} It is easy to check that:
      \[
        \begin{split}
          \langle M_u\mu^\alpha ,\mu^\beta\rangle&=\langle \mu^\alpha,\mu^{\beta+e_u}\rangle\\
          &=c_\alpha \delta_{\alpha, \beta+e_u}\\
          &=c_\alpha \delta_{\alpha-e_u, \beta}\\
&=\frac{c_\alpha}{c_{\alpha-e_u}}c_{\alpha-e_u}\delta_{\alpha-e_u, \beta}\\
&= \frac{c_\alpha}{c_{\alpha-e_u}}\langle\mu^{\alpha-e_u},\mu^\beta\rangle
\end{split}
\]

This concludes the proof.
      \end{proof}

\section{Arveson space and Schur multipliers}
\setcounter{equation}{0}
\label{Vq-Arveson}

For properties of the classical Arveson space, as well as a motivation for its definition, the reader can turn to~\cite{arveson-acta,MR80c:47010}. We now turn to the definition of an Arveson space in our case, which we will denote by $\mathfrak A$, in the present setting.
We start with the definition of the counterpart of the {\em unit ball} in our case:

 \begin{equation}
\Omega_A=\Omega(\mathfrak A)=\left\{x\in\mathbb{H}^*,\,\,\, \sum_{\alpha\in\mathbb N_0^3}|\mu^\alpha(x)|^2\frac{|\alpha|!}{\alpha!}<\infty\right\}.
  \end{equation}

\begin{definition}
The $V_q-$Arveson space is the Hilbert space of absolutely convergent $V_q$-series on the unit ball
  \begin{equation}
\sum_{\alpha\in\mathbb N_0^3}\mu^\alpha(x)c_\alpha
\end{equation}
with
\begin{equation}
\langle f,f\rangle=\sum_{\alpha\in\mathbb N_0^3}\frac{\alpha!}{|\alpha|!}|c_\alpha|^2<\infty,
  \end{equation}
  i.e.:
  \begin{equation}
  \mathfrak A=\{ f(x)=\sum_{\alpha\in\mathbb N_0^3}\mu^\alpha(x)c_\alpha | \, \,f {\rm \,abs.\, conv. \, with\, } \langle f,f\rangle=\sum_{\alpha\in\mathbb N_0^3}\frac{\alpha!}{|\alpha|!}|c_\alpha|^2<\infty\}.
  \end{equation}
\end{definition}
Here the notations are the ones of Definition~\ref{mu-alpha}.

\begin{definition}
\label{Averson_Vq_rep_ker}
The reproducing kernel for this Arveson space is:
  \begin{equation}
K_{\mathfrak A}(x,y)=\sum_{\alpha\in\mathbb N_0^3} \, \mu^\alpha(x)\overline{ \mu^\alpha(y)} \,\frac{|\alpha|!}{\alpha!}.
    \end{equation}
\end{definition}
The domain of this kernel is defined by the following:
\begin{proposition}
  Elements of the Arveson space are $V_q$-regular in
  \begin{equation}
  \Omega_A=\Omega(\mathfrak A)=\left\{x\in\mathbb H^*,\,\,\, \sum_{\alpha\in\mathbb N_0^3}|\mu^\alpha(x)|^2\frac{|\alpha|!}{\alpha!}<\infty\right\}.
  \end{equation}
\end{proposition}
\begin{proof}
This follows from the reproducing kernel property, using the same type of arguments as in Section~\ref{Vq-Kernel-Omega}.
  \end{proof}

Evaluation at the origin does not exist, however, we can now define its counterpart by setting, for $f\in  \mathfrak A$ :
\[
C(f)=c_{0,0,0}.
\]

  \begin{theorem}
    The Arveson space is the unique reproducing kernel Hilbert space (up to a multiplicative constant) of $V_q$ power series
    such that
  $M_1,M_2$ and $M_3$ are bounded operators there and satisfy
  \begin{equation}
    \label{structural-equation}
I-\sum_{u=1}^3M_uM_u^*=CC^*.
\end{equation}
Furthermore,
  \begin{equation}
    M_u^*=B_u,\quad u=1,2,3,
  \end{equation}
  in the Arveson space.
  \end{theorem}

  \begin{proof}
We first assume that \eqref{structural-equation} (sometimes called the structural identity) is in force and compute the associated inner product. We proceed in a number of
steps.\\

STEP 1: {\sl We have that
\[
  C^*1=1_{\mathbb H},
\]
(the constant function equal identically to $1$).}\smallskip

Indeed,
\[
\begin{split}
  (C^*1)(y)&=\langle C^*1 , K_{\mathfrak A}(\cdot, y)\rangle_{\mathfrak A}\\
  &=\langle 1_{\mathbb H} , C K_{\mathfrak A}(\cdot, y)\rangle_{\mathbb H}\\
  &=1.
\end{split}
\]
Here we used the fact that $C K_{\mathfrak A}(\cdot, y)=1_{\mathbb H}$.\\

STEP 2: {\sl $M_u^*1=0$}\smallskip

This is \eqref{alpha-u} with $\alpha=(0,0,0)$.\\

STEP 3: {\sl Let $\alpha$ and $\beta$ different from $(0,0,0)$, and assume that \eqref{structural-equation} is in force.
We have
\begin{equation}
  \label{cdert}
1=\sum_{u=1}^3\frac{c_\alpha}{c_{\alpha-e_u}},\quad \alpha\in\mathbb N^3.
\end{equation}}

Indeed,
\[
    \langle\mu^\alpha,\mu^\beta\rangle=\sum_{u=1}^3\langle M_u^*\mu^\alpha,M_u^*\mu^\beta\rangle
  \]
  that is:
  \[
c_\alpha\delta_{\alpha,\beta}=\sum_{u=1}^3\frac{c_\alpha c_\beta}{c_{\alpha-e_u}c_{\beta-e_u}}\langle\mu^{\alpha-e_u}\mu^{\beta-e_u}\rangle
\]
or, equivalently \eqref{cdert} holds.\\

STEP 4: {\sl The inner product is that of the Arveson space.}\smallskip

From \eqref{cdert} we have that the sequence $\displaystyle d_\alpha=\frac{1}{c_\alpha}$ satisfies
\[
d_\alpha=\sum_{u=1}^3d_{\alpha-e_u},\quad \alpha\in\mathbb N^3.
\]
Together with $d_{0,0,0}=1$ we get $\displaystyle c_\alpha=c_{0,0,0}\frac{\alpha!}{|\alpha|!}$.\\

The converse is proven by reading these arguments backwards, with $\displaystyle c_\alpha=\frac{\alpha!}{|\alpha|!}$.

  \end{proof}

  \begin{proposition}
Assume that $S$ is a $\mathbb H^{n\times m}$-valued multiplier defined on $\Omega_A$. Then,
  \begin{equation}
    \label{waterloo}
    \left(M_S^*(K_{\mathfrak A}(\cdot, y))\xi\right)(x)=\sum_{\alpha\in\mathbb N_0^3}\frac{|\alpha|!}{\alpha!}\mu^\alpha(x)\left(S\star_{V_q}\mu^\alpha)(y)
    \right)^*\xi,\quad \forall \xi\in\mathbb H^n.    
    \end{equation}
    \end{proposition}

\begin{proof}
  Since $M_S^*$ is continuous we obtain: 
  \[
  \begin{split}
    \langle M_S^*(K_{\mathfrak A}(\cdot, y)\xi),K_{\mathfrak A}(\cdot, x)\eta\rangle_{\mathfrak A}&=\langle K_{\mathfrak A}(\cdot, y)\xi,S\star_{V_q}
    K_{\mathfrak A}(\cdot ,x)\eta\rangle_{\mathfrak A}\\
    &=\sum_{\alpha\in\mathbb N_0^3}\frac{|\alpha|!}{\alpha!}\mu^\alpha(x)\langle K_{\mathfrak A}(\cdot, y)\xi, S\eta\star_{V_q}\mu^\alpha\rangle_{\mathfrak A}\\
    &=\sum_{\alpha\in\mathbb N_0^3}\frac{|\alpha|!}{\alpha!}\mu^\alpha(x)\overline{\xi^*\, S\eta\star_{V_q}\mu^\alpha(y)}\\
    &=\eta^*\sum_{\alpha\in\mathbb N_0^3}\frac{|\alpha|!}{\alpha!}\mu^\alpha(x)\left(S\, \xi\star_{V_q}\mu^\alpha(y)\right)^*\\
    &=\eta^*\sum_{\alpha\in\mathbb N_0^3}\frac{|\alpha|!}{\alpha!}\mu^\alpha(x)\left(S\star_{V_q}\mu^\alpha)(y)\right)^*\xi.
  \end{split}
  \]
  \end{proof}

\begin{theorem}
  A $\mathbb H^{n\times m}$-valued function $S$ is a contractive multiplier from $\mathfrak A^m$ into $\mathfrak A^n$ if and
  only if the $\mathbb H^{n\times n}$-valued kernel
  \begin{equation}
    \label{paris}
    K(x,y)=
    \sum_{\alpha\in\mathbb N_0^3}\frac{|\alpha|!}{\alpha!}\left\{\mu^\alpha(x)\overline{\mu^\alpha(y)}I_n-
\left(\mu^\alpha\star_{V_q} S\right)(x)\left(\mu^\alpha\star_{V_q} S\right)(y)^*\right\}
\end{equation}
is positive definite in $\Omega_A=\Omega(\mathfrak A)$.

\end{theorem}

\begin{proof} To simplify notation we assume $n=m=1$.
  Assume first that $M_S$ is a contraction. Then, $I_{\mathfrak A}-M_SM_S^*$ is a positive operator from $\mathfrak A$ into itself. The positivity of the
  kernel \eqref{paris} follows then from the formula
  \[
  \langle( (I_{\mathfrak A}-M_SM_S^*)K_{\mathfrak A}(\cdot, y),K_{\mathfrak A}(\cdot, x)\rangle_{\mathfrak A}=
  \sum_{\alpha\in\mathbb N_0^3}\frac{|\alpha|!}{\alpha!}\left\{\mu^\alpha(x)\overline{\mu^\alpha(y)}I_n-
\left(\mu^\alpha\star_{V_q} S\right)(x)\left(\mu^\alpha\star_{V_q} S\right)(y)^*\right\}
  \]
  which in turn is obtained from \eqref{waterloo}.\smallskip

  Conversely, if the kernel \eqref{paris} is positive definite in $\Omega(\mathfrak A)$, the right linear span of the pairs of functions
  \[
(K_{\mathfrak A}(\cdot, y)))(x),\sum_{\alpha\in\mathbb N_0^3} \frac{|\alpha|!}{\alpha!}(\mu^\alpha(x)\left(S\star_{V_q}\mu^\alpha)(y)\right)^*)\in \mathfrak A\times \mathfrak A
  \]
  defines a densely defined contractive relation, which extends therefore to the graph of an everywhere defined contraction; the adjoint of this contraction
  is $M_S$.
\end{proof}

In view of the structural identity \eqref{cdert} the Fueter-like variables are Schur multipliers. In the next section we present another important Schur multiplier, which is the
counterpart here of an elementary Baschke factor.

\section{$V_q$ Blaschke-type factors}
\setcounter{equation}{0}
\label{Vq-Blaschke}

We denote by $\Omega_1(\mathfrak A)$ the set of $x\in\mathbb H^*$ such that
\[
\sum_{u=1}^3|\mu_u(x)|^2<1.
\]
By inequalities  \eqref{ineq987} we have
\[
\left\{x\in\mathbb H^*\,;\, 3x_0^2+x_1^2+x_2^2+x_3^2<1\right\}\subset \Omega_1(\mathfrak K).
  \]
\begin{remark}
It is essential to consider elements of $\mathbb H^*$, for example the element $x=(\frac{1}{3},0,0,0) \in \left\{x\in\mathbb R^4\,;\, 3x_0^2+x_1^2+x_2^2+x_3^2<1\right\}$. However, $\sum_{u=1}^3|\mu_u(x)|^2$ is not defined. 
  
	It is also interesting to note that the above ellipsoid also appears in \cite{MR2124899}.
\end{remark}
  \begin{theorem}
    Let $a\in\Omega_1(\mathfrak A)$ and set $\mu(a)=\begin{pmatrix}\mu_1(a)&\mu_2(a)&\mu_3(a)\end{pmatrix}$. Then
    the multiplication operator by $\mu(a)$ on the left is a strict contraction 
    from $\mathfrak A^3$ into $\mathfrak A$ and the map
    \begin{equation}
      B_a(x)=(1-\mu(a)\mu(a)^*)^{1/2}
      (1-\mu(x)\mu(a)^*)^{-\star_{V_q}}\star_{V_q}(\mu(x)-\mu(a))(I_3-\mu(a)^*\mu(a))^{-1/2}
\label{blabla}
\end{equation}
     is a Schur multiplier from $\mathfrak A^3$ into $\mathfrak A$.
   \end{theorem}

   \begin{proof}
     The proof follows the proofs in \cite[Proposition 4.1 p. 11]{at-ieot} and \cite[Theorem 4.7 p. 146]{MR2124899},
     and is briefly outlined. In the complex setting we also refer to \cite{rudin-ball}, where a different,
     but equivalent expression
     is given for $B_a$. We set
     \[
J=\begin{pmatrix} I_{\mathfrak A}&0\\0&-I_{{\mathfrak A}^3}\end{pmatrix}.
       \]

       STEP 1: {\sl The operator-matrix
         \[
           H(a)=\begin{pmatrix}(I_{\mathfrak A}-M_{\mu(a)}M_{\mu(a)}^*)^{-1/2}&-M_{\mu(a)}(I_{{\mathfrak A}^3}-M_{\mu(a)}^*M_{\mu(a)})^{-1/2}\\
               -M_{\mu(a)}^*(I_{{\mathfrak A}}-M_{\mu(a)}M_{\mu(a)}^*)^{-1/2}&           (I_{{\mathfrak A}^3}-M_{\mu(a)}^*M_{\mu(a)})^{-1/2}
               \end{pmatrix}
           \]
           is $J$-unitary, i.e.
           \[
H(a)JH(a)^*=H(a)^*JH(a)=J.
\]}\smallskip

See \cite{Dym_CBMS} for the general case where $M_{\mu(a)}$ is replaced by an arbitrary strict contraction between Hilbert spaces.
The operator-matrix $H(a)$ is called the Halmos extension of $M_{\mu(a)}$.\\

STEP 2: {\sl The operator
  \[
\begin{pmatrix}I_{\mathfrak A}& M_{\mu}\end{pmatrix}H(a)JH(a)^*\begin{pmatrix}I_{\mathfrak A}\\ M_{\mu}^*\end{pmatrix}
\]
is non-negative.}\smallskip

Indeed, $H(a)JH(a)^*=J$ by the previous step, and $I_{\mathfrak A}-M_{\mu}M_{\mu}^*\ge$ by \eqref{cdert}.\\

STEP 3: {\sl $B_a$ is a Schur multiplier.}\smallskip

It suffices to write
\[
  \begin{split}
    \begin{pmatrix}I_{\mathfrak A}& M_{\mu}\end{pmatrix}H(a)&=\\
    &\hspace{-2cm}=\begin{pmatrix}
      (I_{\mathfrak A}-M_{\mu}M_{\mu(a)}^*)(I_{\mathfrak A}-M_{\mu(a)}M_{\mu(a)}^*)^{-1/2}&
      (M_{\mu}-M_{\mu(a)})(I_{{\mathfrak A}^3}-M_{\mu(a)}^*M_{\mu(a)})^{-1/2}\end{pmatrix}\\
    &\hspace{-2cm}=(\sqrt{1-\mu(a)^*\mu(a)})^{1/2}
    (I_{\mathfrak A}-M_{\mu}M_{\mu(a)}^*)\begin{pmatrix} I_{\mathfrak A}&M_{B_a}\end{pmatrix}.
    \end{split}
  \]
  This concludes the proof.
         \end{proof}
\section{$V_q$ Rational Functions} 
\setcounter{equation}{0}
\label{V_q-Rational}
The notion of rational function is important in Schur analysis, and leads to state space representations of linear systems.
The study of hypercomplex rational functions using the Cauchy-Kowaleski extension theorem originates with the work of Laville;
see \cite{Laville}.
\begin{definition} 
The $\mathbb H^{n\times m}$-valued $V_q$-regular function on an open domain $\Omega\subset\mathbb{H}^*$ is called $V_q$-rational if its restriction to $x_0=0$
 can be written as
 \begin{equation}
   \label{r010203}
R(0,x_1,x_2,x_3)=D+C(I_N-\sum_{u=1}^3x_uA_u)^{-1}(\sum_{u=1}^3x_uB_u)
\end{equation}
where $D\in\mathbb H^{n\times m}$, $C\in \mathbb H^{n\times N}$, $A_1,A_2,A_3\in\mathbb H^{N\times N}$ and $B_1,\ldots, B_3
\in\mathbb H^{N\times m}$.

\end{definition}

Equivalenty, taking the $V_q$-extension, we can write
\begin{equation}
  R(x)  =D+C(I-\mu(x)A)^{-\star_{V_q}}\star_{V_q}\mu(x)B
\end{equation}
where
\[
  A=\begin{pmatrix}A_1\\ A_2\\ A_3\end{pmatrix}\quad{\rm and}\quad
  B=\begin{pmatrix}B_1\\ B_2\\ B_3\end{pmatrix}.
\]

\begin{proposition}
  The Blaschke factor \eqref{blabla} is a $V_q$-rational function on $\Omega_1(\mathfrak A)$.
\end{proposition}

\begin{proof} We follow the proof of \cite[Proposition 4.1.p. 12]{at-ieot} and first recall that for a contraction
  $K\in\mathbb H^{s\times r}$ it holds that
    \[
      K(I_r-K^*K)^{1/2}=(I_s-KK^*)^{1/2}K\quad {\rm and}\quad      K^*(I_s-KK^*)^{1/2}=(I_r-K^*K)^{1/2}K^*.
    \]
    These equalities are used in the computations below. \smallskip
  
  To see that $B_a$is $V_q$-rational we write (with $\mu(x)=\begin{pmatrix}x_1&x_2&x_3\end{pmatrix}$)
  \[
    \begin{split}
      B_a(0,x_1,x_2,x_3)&=\\
      &\hspace{-2cm}=B_a(0,0,0,0)+B_a(0,x_1,x_2,x_3)-B_a(0,0,0,0)\\
      &\hspace{-2cm}=-\underbrace{(1-\mu(a)\mu(a)^*)^{1/2}
      \mu(a)(I_3-\mu(a)^*\mu(a))^{-1/2}}_{\mu(a)}+\\
    &\hspace{-1.5cm}+(1-\mu(a)\mu(a)^*)^{1/2}(1-\mu(x)\mu(a)^*)^{-1}
    (\mu(x)-\mu(a))(I_3-\mu(a)^*\mu(a))^{-1/2}+\\
      &\hspace{-1.5cm}      +(1-\mu(a)\mu(a)^*)^{1/2}
      \mu(a)(I_3-\mu(a)^*\mu(a))^{-1/2}\\
      &\hspace{-2cm}=-\mu(a)+\\
      &\hspace{-1.5cm}+(1-\mu(a)\mu(a)^*)^{1/2}(1-\mu(x)\mu(a)^*)^{-1}\times\\
      &\hspace{-1.cm}\times
\left\{\mu(x)-\mu(a)+(1-\mu(x)\mu(a)^*)\mu(a)\right\}(I_3-\mu(a)^*\mu(a))^{-1/2}\\
      &\hspace{-2cm}=-\mu(a)+(1-\mu(a)\mu(a)^*)^{1/2}(1-\mu(x)\mu(a)^*)^{-1}\mu(x)(I_3-\mu(a)^*\mu(a))^{1/2}.
      \end{split}
    \]
    which is of the form \eqref{r010203} with
    \begin{equation}
      \label{T-iso}
   T=   \begin{pmatrix}A_1&B_1\\
        A_2&B_2\\
        A_3&B_3\\
        C&D\end{pmatrix}=\begin{pmatrix}\mu(a)^*&(I_3-\mu(a)^*\mu(a))^{1/2}\\
        (1-\mu(a)\mu(a)^*)^{1/2}&-\mu(a)\end{pmatrix}.
      \end{equation}
    \end{proof}

    \begin{remark}
      {\rm We note that \eqref{T-iso} is co-isometric. Existence of a co-isometric realization is a general property of Schur
        multipliers, and will be considered elsewhere. Here, $T$ is in fact unitary.}
      \end{remark}

Rational functions were defined in the setting of the Cauchy operator in \cite{MR2124899}. When restricted to $x_0=0$, one obtains
the same class of functions, namely functions of the of form \eqref{r010203}. $V_q$-rational functions are stable under sum,
$V_q$-product and $V_q$-inversion when sizes are compatible. The arguments are the same as in the paper \cite{MR2124899},
to which we refer the reader.\\

\section{Conclusions and Future Endeavors}
\setcounter{equation}{0}
We are now extending this work to analyze a theory of Schur-Agler functions (see \cite{agler-hellinger} and see \cite{MR2240272} for the
Fueter case) and Schur multipliers as well as
expand the rational function theory in this context.\smallskip

Consider a system of linear ordinary differential equations with variables coefficients of the form
\[
  \begin{split}
V_j(f)=\frac{\partial f_j}{\partial x_0}&-\sum_{u=1}^3\sum_{k=0}^3 a_{k,u}(x)\frac{\partial f_k}{\partial x_u}=0,\quad j=0,1,2,3,
  \end{split}
  \]
  where the $a_{k,u}$ are real analytic on some open subset of the real line, and the system of equations is denoted by $V$. Thanks to the Cauchy-Kovalevskaia theorem one can define
  $V$-Fueter variables as the CK extensions of the functions $x_u$, $u=1,2,3$ and monomials as CK-extensions of the monomials
  $x^\alpha$. It seems difficult though to develop a general theory which goes beyond some trivial facts. The cases considered
  in \cite{MR2124899} (for the Cauchy-Fueter operator, factorizing the $\mathbb R^4$ Laplacian) and in \cite{MR3795166} 
   associated to the operator
  \[
    \frac{\partial^2}{\partial x_0^2}-\frac{\partial^2}{\partial x_1^2}-\frac{\partial^2}{\partial x_2^2}+
    \frac{\partial^2}{\partial x_3^2}
    \]
    studied in the setting of split quaternion, together with the present analysis for the $V_q$-operator exhibit how different
    each specific case can be.\\

    \bibliographystyle{plain}

\end{document}